\documentclass{amsart}

\RequirePackage{fix-cm}
\usepackage{graphicx}
\usepackage{amsmath, amsopn,amstext,amscd,amsfonts,amssymb}
\usepackage{dsfont}
\usepackage{comment}
\usepackage[active]{srcltx}
\usepackage{graphicx, epsfig, subfig}

\usepackage{textcomp}

\usepackage{algorithm}

\makeatletter
\def\BState{\State\hskip-\ALG@thistlm}
\makeatother

\def\downbar#1{
\setbox10=\hbox{$#1$}
            \dimen10=\ht10 \advance\dimen10 by 2.5pt
            \ifdim \dimen10<15pt %equals approximately 0.5cm
               \advance\dimen10 by -0.5pt
               \dimen11=\dimen10
               \advance\dimen10 by 2.5pt
               \lower \dimen11
            \else \lower \ht10 \fi
            \hbox {\hskip 1.5pt \vrule height \dimen10 depth \dp10}}
\def\upbar#1{
\setbox10=\hbox{$#1$}
            \dimen10=\ht10 \advance\dimen10 by \dp10 \advance\dimen10 by 2.5pt
            \ifdim \dimen10<15pt %equals approximately 0.5cm
                \advance\dimen10 by 2pt \fi
            \raise 2.5pt \hbox {\hskip -1.5pt \vrule height \dimen10}}

\usepackage{multicol}
\usepackage{colortbl}
\usepackage{exscale}
\usepackage{amssymb,latexsym,amsthm,amsfonts,color,fancyhdr,mathrsfs}
\usepackage{lscape}
\usepackage{rotating}
\usepackage{caption}

\newtheorem{definition}{\bf Definition}[section]
\newtheorem{theorem}{\bf Theorem}[section]

\newtheorem{lemma}{\bf Lemma }[section]

\newtheorem{remark}{\bf Remark}[section]

\numberwithin{equation}{section}

\begin{document}

\title[On another extension of coherent pairs of measures]{On another extension of coherent pairs of measures}
\author{K. Castillo}
\address{CMUC, Department of Mathematics, University of Coimbra, 3001-501 Coimbra, Portugal}
\email{kenier@mat.uc.pt}
\author{D. Mbouna}
\address{CMUC, Department of Mathematics, University of Coimbra, 3001-501 Coimbra, Portugal}
\email{dmbouna@mat.uc.pt}

\date{\today}
\begin{abstract}
Let $M$ and $N$ be fixed non-negative integer numbers and let
$\pi_N$ be a polynomial of degree $N$. Suppose that $(P_n)_{n\geq0}$ and $(Q_n)_{n\geq0}$
are two orthogonal polynomial sequences such that
%their derivatives of orders $k$ and $m$ (respectively) satisfy the structure relation
$$
\pi_N(x)\,P_{n+m}^{(m)}(x)= \sum_{j=n-M}^{n+N}r_{n,j}Q_{j+k}^{(k)}(x)\quad (n=0,1,\ldots)\,,
$$
where $r_{n,j}$ are complex number independent of $x$.
It is shown that under natural constraints, $(P_n)_{n\geq0}$ and $(Q_n)_{n\geq0}$
are semiclassical orthogonal polynomial sequences.
Moreover, their corresponding moment linear functionals are related by a rational
modification in the distributional sense.
This leads to the concept of $\pi_N-$coherent pair with index $M$ and order $(m,k)$.\end{abstract}
\maketitle
\section{Introduction}
In the framework of the theory of orthogonal polynomials
---for an updated reference on this subject we refer the reader to Ismail's book \cite{I2005}---, the
concept of coherent pair of measures as well as its multiple generalizations have been a subject of increasing research interest along the last decades. This concept was introduced by Iserles et al. \cite{IKNS1991}
motivated by the theory of polynomial approximation with respect to certain Sobolev inner products.
In \cite{JP2013,JMPN2015}, the notion of $(M,N)-$coherent pair, and of $(M,N)-$coherent pair of order $(m,k)$,
were introduced as extensions of most of the concepts of coherence up to that time.
More precisely, given two monic orthogonal polynomial sequences (OPS), $(P_n)_{n\geq0}$ and $(Q_n)_{n\geq0}$,
we say that $\big((P_n)_{n\geq0},(Q_n)_{n\geq0}\big)$ is an $(M,N)-$coherent pair of order $(m,k)$ if there exist two non-negative integer numbers $M$ and $N$, and sequences of complex numbers
$(a_{n,j})_{n\geq0}$ ($j=0,1,\ldots,M$) and $(b_{n,j})_{n\geq0}$ ($j=0,1,\ldots,N$) such that,
under natural assumptions on the coefficients $a_{n,j}$ and $b_{n,j}$, the structure relation
$$
\sum_{j=0}^{M} a_{n,j} P_{n-j}^{[m]} (x) = \sum_{j=0}^{N} b_{n,j} Q_{n-j}^{[k]} (x) \quad (n=0,1,\ldots)
$$
holds. Here and subsequently, we use the notation
$$
%P_{n}^{[m]}(x):=\frac{{\rm d}^m}{{\rm d}x^m}\left(\frac{P_{n+m}(x)}{(n+1)_m}\right)
P_{n}^{[m]}(x):=\frac{1}{(n+1)_m}\frac{{\rm d}^m}{{\rm d}x^m}\,P_{n+m}(x)
$$
($Q_{n}^{[k]}$ is defined in the same way), where for any positive real number $\alpha$,
$(\alpha)_n$ denotes the Pochhammer symbol defined by
$$(\alpha)_0:=1\;,\quad (\alpha)_n:=\alpha(\alpha+1)\cdots(\alpha+n-1)\quad\mbox{\rm if}\quad n\in\mathbb{N}\;.$$
Note that $P_{n}^{[m]}$ is a normalization of the derivative of order $m$ of $P_{n+m}$
defined so that it becomes a monic polynomial of degree $n$.
Let $\textbf{u}$ and $\textbf{v}$ be the moment regular functionals with respect to which
$(P_n)_{n\geq0}$ and $(Q_n)_{n\geq0}$ are orthogonal.
It follows from the results in \cite{P2006,JP2008,JP2013,JMPN2015} that if $m=k$ then $\textbf{u}$ and $\textbf{v}$ are connected by a rational transformation (in the distributional sense), i.e., there exist nonzero polynomials $\Phi$ and $\Psi$ such that $\Phi\textbf{u}=\Psi\textbf{v}$. Otherwise if $m \neq k$ then $\textbf{u}$ and $\textbf{v}$ are still connected by a rational transformation and, in addition, they are semiclassical functionals, i.e., there exist nonzero polynomials $\Phi_1$, $\Psi_1$, $\Phi_2$, and $\Psi_2$ such that
$$ D(\Phi_1{\bf u})=\Psi_1{\bf u}\;,\quad D(\Phi_2{\bf v})=\Psi_2{\bf v}\;.$$
In agreement with the `algebraic theory' of OPS introduced by Maroni \cite{M1991},
the left product of a polynomial $\Phi$ by a moment functional $\textbf{w}$ is
the functional, $\Phi{\bf w}$, defined by $\langle \Phi\textbf{w},p\rangle:=\langle \textbf{w},\Phi p\rangle$ for each polynomial $p$, whereas the derivative of ${\bf w}$, $D{\bf w}$, is defined by $\langle D\textbf{w},p\rangle:=-\langle \textbf{w},p'\rangle$, for each polynomial $p$.
As usual, $\langle\cdot,\cdot\rangle$ means the duality bracket, so that $\langle {\bf w},p\rangle$ is the action of the functional ${\bf w}$ over the polynomial $p$.

In this work we modify the left-hand side of the above structure relation, and consider the following one:
\begin{equation}\label{1b}
\pi_N(x)P_n^ {[m]}(x) = \sum_{j=n-M}^{n+N} c_{n,j} Q_j ^{[k]}(x)\quad (n=0,1,\ldots)\,,
\end{equation}
where $M$ and $N$ are fixed non-negative integer numbers, $\pi_N $ is a monic polynomial of degree $N$
(hence $c_{n,n+N}=1$ for each $n$), and we consider the convention $Q_j\equiv0$ if $j<0$.
Further, we will assume that the following conditions hold:
\begin{equation}\label{cond1}
c_{n,n-M} \neq 0\quad\mbox{\rm if}\quad n\geq M\;.
\end{equation}
Maroni and Sfaxi \cite{MS2000} considered the case $(m,k)=(0,1)$
and called the pair $\big((P_n)_{n\geq0},(Q_n)_{n\geq0}\big)$ fulfilling the structure relation
(\ref{1b}) whenever  $(m,k)=(0,1)$ a
$\pi_N-$coherent pair with index $M$. This motivates the following

\begin{definition}
Let $M$ and $N$ be non-negative integer numbers and let
$\pi_N$ be a monic polynomial of degree $N$. If $(P_n)_{n\geq0}$ and $(Q_n)_{n\geq0}$
are two monic OPS such that their normalized derivatives of orders $m$ and $k$ (respectively)
satisfy $(\ref{1b})$--$(\ref{cond1})$, we call $\big((P_n)_{n\geq0},(Q_n)_{n\geq0}\big)$,
as well as the corresponding pair $({\bf u},{\bf v})$ of regular functionals, a
$\pi_N-$coherent pair with index $M$ and order $(m,k)$.
\end{definition}

Besides \cite{MS2000}, many other instances of the structure relation (\ref{1b}) were considered previously by several authors. For instance, the case $N=0$ (i.e., $\pi_N\equiv1$ and $M$, $m$, and $k$ being arbitrary) fits into the theory of $(M,0)-$coherent pairs of order $(m,k)$, described at the begin of this introduction. Also, whenever $(m,k)=(1,0)$ and $(P_n)_{n\geq0}\equiv (Q_n)_{n\geq0}$, $(\ref{1b})$ becomes a characterization of semiclassical OPS due to Maroni \cite{M1987,M1991}.
Note that for $N\leq2$ and $M=0$, this reduces to the well known Al-Salam-Chihara characterization of the classical OPS \cite{A-SC1972}.
The case $k=0$ ($M$, $N$ and $m$ being arbitrary) was considered by Bonan et al. \cite{BLN1987} in the framework of orthogonality in the positive-definite sense, i.e., whenever the orthogonality of each of the involved OPS is considered with respect to positive Borel measures.
In the special case $m=1$, a complementary approach to the case considered in \cite{BLN1987} was presented in \cite{MBP1993},
in the framework of the so-called regular (or formal) orthogonality.

It is a remarkable fact that in all the previous works the involved OPS and their corresponding regular moment linear functionals are semiclassical. Thus, a major question is to analyze whether the OPS involved in a $\pi_N-$coherent pair with index $M$ and order $(m,k)$ are semiclassical, and in such a case to determine the relations between the corresponding regular moment linear functionals.
This will be treated in Section \ref{Section-main}.
As an application, in Section \ref{Section-example}, we present an alternative approach
to a recent result due to Griffin \cite{G2016},
which fits into $\pi_1-$coherence with index $1$ and order $(1,0)$.

\section{Main results}\label{Section-main}

In this section we establish the semiclassical character of the OPS and their associated regular functionals
involved in a $\pi_N-$coherent pair with index $M$ and order $(m,k)$.
Our approach is based upon the algebraic theory of orthogonal polynomials
developed by Maroni \cite{M1985,M1991}. We denote by $\mathcal{P}$ the vector space of all (complex) polynomials
and by $\mathcal{P}^*$ its algebraic dual space. $\mathcal{P}$ may be endowed with a topology
(indeed, an appropriate strict inductive limit topology) such that the algebraic and
the topological dual spaces of $\mathcal{P}$ coincide, that is, $\mathcal{P}^*=\mathcal{P}'$.
Given a simple set of polynomials $(R_n)_{n\geq0}$
(meaning that each $R_n\in\mathcal{P}$ and $\deg R_n=n$ for each $n=0,1,\ldots$),
the corresponding dual basis is a sequence of linear functionals
${\bf e}_n:\mathcal{P}\to\mathbb{C}$ such that
$$
\langle{\bf e}_n,R_j\rangle:=\delta_{n,j}\quad(n,j=0,1,\ldots)\;,
$$
where $\delta_{n,j}$ denotes the Kronecker's symbol.
In particular, if $(R_n)_{n\geq0}$ is a monic OPS with respect to ${\bf w}\in\mathcal{P}'$, i.e.,
there exists a sequence of nonzero complex numbers $(k_n)_{n\geq0}$ such that the orthogonality conditions
$$
\langle{\bf w},R_jR_n\rangle:=k_n\delta_{j,n}\quad(j,n=0,1,\ldots)
$$
hold, then the corresponding dual basis is explicitly given by
$$
{\bf e}_n=k_n^{-1}R_n{\rm w}\quad (n=0,1,\ldots)\;.
$$

\begin{lemma}\label{lemma-main}
Let $\big((P_n)_{n\geq0},(Q_n)_{n\geq0}\big)$ be a $\pi_N-$coherent pair with index $M$ and order $(m,k)$,
so that $(\ref{1b})$--$(\ref{cond1})$ hold. Set
\begin{align}
& \psi(x;n):=\sum_{j=n-N}^{n+M}\frac{(-1)^m(j+1)_m\; c_{j,n}}{\langle {\bf u}, P_{m+j} ^2\rangle} P_{m+j}(x)\;, \label{si} \\
& \phi(x;n,j):= \frac{(-1)^k(n+1)_k}{\langle{\bf v},Q_{n+k}^2\rangle}\sum_{\ell=0}^{N-j}\binom{k+N}{\ell}
\binom{N-\ell}{N-j-\ell} \pi_N^{(\ell)}(x) Q_{n+k}^{(N-j-\ell)}(x)\,, \label{fi}
\end{align}
for all $\,n=0,1,\ldots$ and $j=0,1,\ldots,N$, so that
$$\deg\psi(\cdot;n)=m+n+M\;,\quad\deg\phi(\cdot;n,j)=k+n+j\;.$$
Let ${\bf u}$ and ${\bf v}$ be the regular functionals with respect to which
$(P_n)_{n\geq0}$ and $(Q_n)_{n\geq0}$ are orthogonal.
Then the following functional equations hold:
\begin{align}
& \psi(\cdot;n){\bf u}=D^{m-k-N}\left(\sum_{j=0}^{N}\phi(\cdot;n,j)D^{j}{\bf v}\right)
& \mbox{\rm if}\quad m \geq k+N\;, \label{2a} \\
& D^{k+N-m}\big(\psi(\cdot;n){\bf u}\big)=\sum_{j=0}^{N}\phi(\cdot;n,j)D^j{\bf v}
& \mbox{\rm if}\quad m< k+N\;, \label{2b}
\end{align}
for all $n=0,1,\ldots$.
\end{lemma}

\begin{proof}
Let $({\bf a}_n)_{n\geq0}$, $({\bf b}_n)_{n\geq0}$, $({\bf a}_n^{[m]})_{n\geq0}$,
and $({\bf b}_n^{[k]})_{n\geq0}$ be the dual basis corresponding to the simple sets of polynomials
$(P_n)_{n\geq0}$, $(Q_n)_{n\geq0}$, $(P_n^{[m]})_{n\geq0}$ and $(Q_n^{[k]})_{n\geq0}$, respectively.
Then
$$
\pi_N {\bf b}_n^{[k]} = \sum_{j=0}^{+\infty} \big\langle \pi_N {\bf b}_n^{[k]}, P_j^{[m]} \big\rangle\, {\bf a}_j^{[m]}
\quad(n=0,1,\ldots)$$
(in the sense of the weak dual topology in $\mathcal{P}'$).
From $(\ref{1b})$, we have
\begin{align*}
\big\langle \pi_N {\bf b}_n^{[k]}, P_j ^{[m]} \big\rangle &= \big\langle b_n ^{[k]}, \pi_N P_j ^{[m]} \big\rangle
=\sum_{\ell=j-M}^{j+N} c_{j,\ell} \big\langle  b_n^{[k]}, Q_\ell^{[k]} \big\rangle \\
&= \begin{cases} c_{j,n} & \mbox{\rm if }\quad n-N \leq j\leq n+M  \\
 0 & \mbox{\rm otherwise\,. } \end{cases}
\end{align*}
Hence
\begin{equation}\label{rr}
\pi_N {\bf b}_n ^{[k]} = \sum_{j=n-N} ^{n+M} c_{j,n} {\bf a}_j ^{[m]} \quad(n =0,1,\ldots)\;.
\end{equation}
Considering the $m$-th derivative on both sides of this equation
and taking into account that
$D^m\big({\bf a}_j ^{[m]}\big)=(-1)^m(j+1)_m{\bf a}_{j+m}$, we obtain
\begin{equation}\label{basis}
D^m \big(\pi_N {\bf b}_n ^{[k]}  \big) = \psi(\cdot;n)\textbf{u} \quad(n=0,1,\ldots)\;,
\end{equation}
where $\psi(\cdot;n)$ is defined by (\ref{si}).
Notice that the condition (\ref{cond1}) ensures that $\deg\psi(\cdot,n)=M+m+n$ for each $n=0,1,\ldots$.
Using the Leibniz rule for the derivative of the left product of a functional by a polynomial,
and taking into account that $\pi_N^{(j)}=0$ if $j>N$, as well as
$$
D^k\big({\bf b}_n ^{[k]}\big)=(-1)^k(n+1)_k{\bf b}_{n+k}
=(-1)^k(n+1)_k\langle{\bf v},Q_{n+k}^2\rangle^{-1}Q_{n+k}{\bf v}\,,
$$
we deduce
\begin{align*}
&D^{k+N} \big(\pi_N {\bf b}_n ^{[k]} \big)\\
&\qquad=\frac{(-1)^k(n+1)_k}{\langle {\bf v}, Q_{n+k}^2 \rangle} \sum_{j=0}^{N} \binom{k+N}{j} \pi_N^{(j)}  D^{N-j} (Q_{n+k} {\bf v}) \\
&\qquad= \frac{(-1)^k (n+1)_k }{\langle {\bf v}, Q_{n+k} ^2 \rangle} \sum_{j=0} ^{N} \sum_{\ell=0} ^{N-j}  \binom{k+N}{j} \binom{N-j}{\ell} \pi_N ^{(j)} Q_{n+k}^{(\ell)} D^{N-j-\ell}{\bf v} \\
&\qquad= \frac{(-1)^k (n+1)_k}{\langle {\bf v}, Q_{n+k} ^2\rangle} \sum_{j=0}^{N} \sum_{\ell=j} ^{N}  \binom{k+N}{j} \binom{N-j}{\ell-j}  \pi_N ^{(j)} Q_{n+k}^{(\ell-j)} D^{N-\ell} {\bf v} \\
&\qquad= \frac{(-1)^k (n+1)_k }{\langle {\bf v}, Q_{n+k}^2\rangle} \sum_{\ell=0}^{N} \sum_{j=0}^{\ell} \binom{k+N}{j} \binom{N-j}{\ell-j} \pi_N^{(j)} Q_{n+k}^{(\ell-j)}  D^{N-\ell}{\bf v}  \\
&\qquad= \sum_{\nu=0}^{N}\left(\frac{(-1)^k (n+1)_k }{\langle {\bf v}, Q_{n+k}^2\rangle}  \sum_{j=0}^{N-\nu} \binom{k+N}{j} \binom{N-j}{N-\nu-j} \pi_N^{(j)} Q_{n+k}^{(N-\nu-j)}\right)  D^{\nu}{\bf v} \;.
\end{align*}
Hence, by (\ref{fi}), we obtain
\begin{equation}\label{EqDkN1}
D^{k+N} \big(\pi_N {\bf b}_n ^{[k]} \big)= \sum_{j=0}^{N}  \phi(\cdot;n,j) D^{j} {\bf v}\;.
\end{equation}
If $m \geq k+N$, we rewrite (\ref{basis}) as
\begin{equation}\label{basis-1}
\psi(\cdot;n){\bf u} =D^{m-k-N } D^{k+N} \big(\pi_N {\bf b}_n ^{[k]}  \big)\quad (n=0,1,\ldots)\;,
\end{equation}
and (\ref{2a}) follows from (\ref{EqDkN1}) and (\ref{basis-1}).
If $m< k+N$, writing
$$
D^{k+N} \big(\pi_N {\bf b}_n ^{[k]} \big)=D^{k-m+N } D^{m} \big(\pi_N {\bf b}_n ^{[k]} \big)
\quad (n=0,1,\ldots)\;,
$$
we see that (\ref{2b}) follows from (\ref{basis}) and (\ref{EqDkN1}).
\end{proof}

Let us first consider the case $m\geq k+N$.

\begin{theorem}\label{theorem-main}
Let $\big((P_n)_{n\geq0},(Q_n)_{n\geq0}\big)$ be a $\pi_N-$coherent pair with index $M$ and order $(m,k)$,
so that $(\ref{1b})$--$(\ref{cond1})$ holds.
Let ${\bf u}$ and ${\bf v}$ be the regular functionals with respect to which
$(P_n)_{n\geq0}$ and $(Q_n)_{n\geq0}$ are orthogonal.
Suppose $m \geq k+N$. Assume further that $m>k$ whenever $N=0$.
For each $i=0,\ldots,m-k$ and $n=0,1,\ldots$, let
\begin{equation}\label{EqVarphi1}
\varphi(x;n,i):=\sum_{\substack{j+\ell=i \\ 0\leq j\leq N \\ 0\leq\ell\leq M}}
\binom{m-k-N}{\ell} \big(\phi(x;n,j)\big)^{(m-k-N-\ell)}\;,
\end{equation}
$\phi(\cdot;n,j)$ being the polynomial introduced in $(\ref{fi})$.
Let $\mathcal{A}(x)$ be the polynomial matrix of order $m-k+1$ defined by
$$
\mathcal{A}(x):=\big[\varphi(x;n,j)\big]_{n,j=0}^{m-k}\;.
$$
Let $\mathcal{A}_1(x)$ (resp., $\mathcal{A}_2(x)$) be the matrix obtained
by replacing the first (resp., the second) column of
$\mathcal{A}(x)$ by $\big[\psi(x;0),\psi(x;1),\cdots,\psi(x;m-k)]^t$,
and set
$$
A(x):=\det\mathcal{A}(x)\;,\quad A_1(x):=\det\mathcal{A}_1(x)\;,\quad A_2(x):=\det\mathcal{A}_2(x)\;.
$$
Assume that the polynomial $A(x)$ does not vanishes identically.
Then
\begin{equation}\label{EqAvA1u}
A{\bf v}=A_1{\bf u}\;,\quad AD{\bf v}=A_2{\bf u}\;,
\end{equation}
hence ${\bf u}$ and ${\bf v}$ are semiclassical functionals related by a rational transformation.
Moreover, ${\bf u}$ and ${\bf v}$ fulfill the following equations:
\begin{equation}\label{FucEquv}
D(AA_1{\bf u})=\big(2A'A_1+AA_2\big){\bf u}\;,\quad D(AA_1{\bf v})=\big((AA_1)'+AA_2\big){\bf v}\;.
\end{equation}
\end{theorem}

\begin{proof}
By (\ref{2a}) and Leibniz rule, we have
$$
\psi(\cdot;n){\bf u}
=\sum_{j=0}^{N}\sum_{\ell=0}^{m-k-N} \binom{m-k-N}{\ell} \big(\phi(\cdot;n,j)\big)^{(m-k-N-\ell)} D^{j+\ell} {\bf v}\;.
$$
This may be rewritten as
\begin{equation}\label{varphi-uv}
\psi(\cdot;n){\bf u}=\sum_{i=0}^{m-k}\varphi(\cdot;n,i)D^{i}{\bf v}
\quad(n=0,1,\dots)\;,
\end{equation}
where $\varphi(\cdot;n,i)$ is the polynomial introduced in (\ref{EqVarphi1}).
Taking $n=0,1,\ldots,m-k$ in (\ref{varphi-uv}) we obtain a system with $m-k+1$ equations that can be written as
$$
\left(
\begin{array}{c}
\psi(x;0){\bf u} \\
\psi(x;1){\bf u} \\
\vdots \\
\psi(x;m-k){\bf u}
\end{array}
\right)
=\mathcal{A}(x)
\left(
\begin{array}{c}
{\bf v} \\
D{\bf v} \\
\vdots \\
D^{m-k}{\bf v}
\end{array}
\right)\;.
$$
Solving for ${\bf v}$ and $D{\bf v}$ we obtain (\ref{EqAvA1u}).
Finally, (\ref{FucEquv}) follows from (\ref{EqAvA1u}). %conclude the statement of the theorem.
\end{proof}

\begin{remark}
If $m=k$ and $N=0$, then ${\bf u}$ and ${\bf v}$ are still related by a rational transformation,
but we cannot ensure that they are semiclassical (see \cite{JP2008,JMPN2015}).
\end{remark}

Now, we consider the case $m<k+N$.

\begin{theorem}\label{theorem-mainB}
Let $\big((P_n)_{n\geq0},(Q_n)_{n\geq0}\big)$ be a $\pi_N-$coherent pair with index $M$ and order $(m,k)$,
so that $(\ref{1b})$--$(\ref{cond1})$ holds.
Let ${\bf u}$ and ${\bf v}$ be the regular functionals with respect to which
$(P_n)_{n\geq0}$ and $(Q_n)_{n\geq0}$ are orthogonal.
Assume further that $m < k+N$.
For each $j=0,\ldots,k-m+N$ and $n=0,1,\ldots$, set
\begin{equation}\label{Eqxi1}
\xi(x;n,j):=\binom{k-m+N}{j}\big(\psi(x;n)\big)^{(k-m+N-j)}\;,
\end{equation}
$\psi(\cdot;n)$ being the polynomial introduced in $(\ref{si})$.
Let $\mathcal{B}(x):= \big[b_{i,j}(x)\big]_{i,j=0}^{k-m+2N}$ be the polynomial matrix of order $k-m+2N+1$ defined by
$$
b_{i,j}(x):=\left\{
\begin{array}{ccl}
\phi(x;i,j) & \mbox{\rm if} & 0\leq j\leq N\;, \\ [0.5em]
-\xi(x;i,j-N) & \mbox{\rm if} & N+1\leq j\leq k-m+2N\;,
\end{array}
\right.
$$
$\phi(\cdot;i,j)$ being the polynomial given by $(\ref{fi})$.
Let $\mathcal{B}_1(x)$ (resp., $\mathcal{B}_2(x)$ and $\mathcal{B}_{N+2}(x)$) be the matrix obtained
by replacing the first (resp., the second and the $(N+2)$-th) column of
$\mathcal{B}(x)$ by $\big[\xi(x;0,0),\xi(x;1,0),\cdots,\xi(x;m-k+2N,0)]^t$,
and set
$$
B(x):=\det\mathcal{B}(x)\;,\quad B_j(x):=\det\mathcal{B}_j(x)\;,\quad j\in\{1,2,N+2\}\;.
$$
Assume that the polynomial $B(x)$ does not vanishes identically.
Then
\begin{equation}\label{EqAvB1u}
B{\bf v}=B_1{\bf u}\;,\quad BD{\bf v}=B_2{\bf u}\;,\quad BD{\bf u}=B_{N+2}{\bf u}\;,
\end{equation}
hence ${\bf u}$ and ${\bf v}$ are semiclassical functionals related by a rational transformation.
Moreover, ${\bf u}$ and ${\bf v}$ fulfill the following equations:
\begin{equation}\label{FucEquvB}
D(B{\bf u})=\big(B'+B_{N+2}\big){\bf u}\;,\quad D(BB_1{\bf v})=\big((BB_1)'+BB_2\big){\bf v}\;.
\end{equation}
\end{theorem}

\begin{proof}
By the Leibniz rule, we can rewrite (\ref{2a}) as
$$
\sum_{j=0}^{k-m+N}\xi(\cdot;n,j)D^j{\bf u}=\sum_{j=0}^{N}\phi(\cdot;n,j)D^j{\bf v}
\quad (n=0,1,\ldots)\;.
$$
Taking $n=0,1,\ldots,k-m+2N$, we obtain the following system of $k-m+2N+1$ equations:
$$
\left(
\begin{array}{c}
\xi(x;0,0){\bf u} \\
\xi(x;1,0){\bf u} \\
\vdots \\
\xi(x;k-m+N,0){\bf u} \\
\xi(x;k-m+N+1,0){\bf u} \\
\vdots \\
\xi(x;k-m+2N,0){\bf u}
\end{array}
\right)
=\mathcal{B}(x)
\left(
\begin{array}{c}
{\bf v} \\
D{\bf v} \\
\vdots \\ [0.25em]
D^N{\bf v} \\
D{\bf u} \\
\vdots \\ [0.25em]
D^{k-m+N}{\bf u}
\end{array}
\right)\;.
$$
The theorem follows by solving this system for ${\bf v}$, $D{\bf v}$, and $D{\bf u}$.
\end{proof}

In the case $k=0$ we may state a finer result.
Recall that if ${\bf u}\in\mathcal{P}'$ is a semiclassical functional then the class of ${\bf u}$,
denoted by $\mathfrak{s}_{\bf u}$, is the unique non-negative integer number defined by
$$
\mathfrak{s}({\bf u}):=\min_{(\Phi,\Psi)\in\mathcal{A}_{\bf u}}\max\big\{\deg\Phi-2,\deg\Psi-1\big\}\;,
$$
where $\mathcal{A}_{\bf u}$ is the set of all pairs of nonzero polynomials $(\Phi,\Psi)$ fulfilling
the functional equation $D(\Phi{\bf u})=\Psi{\bf u}$.

\begin{theorem}\label{theorem-main-kzero}
Let $\big((P_n)_{n\geq0},(Q_n)_{n\geq0}\big)$ be a $\pi_N-$coherent pair with index $M$ and order $(m,0)$,
so that the structure relation
$$
\pi_N(x)P_n^ {[m]}(x) = \sum_{j=n-M}^{n+N} c_{n,j} Q_j(x)\quad (n=0,1,\ldots)
$$
holds, where $M$ and $N$ are fixed non-negative integer numbers, $\pi_N $ is a monic polynomial of degree $N$,
and $c_{n,n-M} \neq 0$ if $n\geq M$.
Assume further that $m\geq1$ if $N=0$.
Let ${\bf u}$ and ${\bf v}$ be the regular functionals with respect to which
$(P_n)_{n\geq0}$ and $(Q_n)_{n\geq0}$ are (respectively) orthogonal.
Then ${\bf u}$ and ${\bf v}$ are semiclassical functionals
related by a rational transformation. More precisely, setting
\begin{equation}\label{EqPhikzero1}
\Phi(x;j):=\frac{\langle{\bf v},Q_j^2\rangle\psi(x;j)-\sum_{\ell=0}^{j-1}
\binom{m}{\ell}Q_j^{(\ell)}(x)\Phi(x;\ell)}{j!\binom{m}{j}}\quad(j=0,1,\ldots,m)\;,
\end{equation}
$\psi(\cdot;j)$ being the polynomial introduced in $(\ref{si})$,
then $\deg\Phi(\cdot;0)=M+m$, $\deg\Phi(\cdot;j)\leq M+m+j$ for each $j=1,\ldots,m$,
and the following holds:
\begin{align}
& D\big(\Phi(\cdot;1){\bf u}\big)=\Phi(\cdot;0){\bf u} \label{EqFuc1}\\
& \pi_N{\bf v}=\Phi(\cdot;m){\bf u} \label{EqFuc2} \\
& D\big(\Phi(\cdot;m)\pi_N{\bf v}\big)=\big(\Phi(\cdot;m)'+\Phi(\cdot;m-1)\big)\pi_N{\bf v} \;. \label{EqFuc3}
\end{align}
Moreover, $\mathfrak{s}({\bf u})\leq M+m-1$ and $\;\mathfrak{s}({\bf v})\leq N+M+2(m-1)$.
\end{theorem}

\begin{proof}
Since $k=0$ then ${\bf b}_n^{[k]}\equiv{\bf b}_n^{[0]}={\bf b}_n=\langle {\bf v},Q_n^2\rangle^{-1}Q_n{\bf v}$ for each $n=0,1,\ldots$,
hence relation (\ref{basis}) may be rewritten as
\begin{equation}\label{basis-kzero1}
D^m \big(Q_n \pi_N {\bf v}\big) = \langle {\bf v},Q_n^2\rangle\psi(\cdot;n)\textbf{u} \quad(n=0,1,\ldots)\;,
\end{equation}
where $\psi(\cdot;n)$ is defined by (\ref{si}).
Taking $n=0$, we obtain
\begin{equation}\label{basis-kzero2}
D^m \big(\pi_N {\bf v}\big) = \Phi(\cdot;0)\textbf{u} \;.
\end{equation}
Taking $n=1$ in (\ref{basis-kzero1}) and then applying the Leibniz rule, we deduce
$$
\langle {\bf v},Q_1^2\rangle\psi(\cdot;1)\textbf{u}=
D^m \big(Q_1 \pi_N {\bf v}\big)=mD^{m-1} \big(\pi_N {\bf v}\big)+Q_1D^{m} \big(\pi_N {\bf v}\big)\;.
$$
Hence, by (\ref{basis-kzero2}), we have
\begin{equation}\label{basis-kzero3}
D^{m-1} \big(\pi_N {\bf v}\big) = \Phi(\cdot;1)\textbf{u} \;.
\end{equation}
Thus (\ref{EqFuc1}) follows from (\ref{basis-kzero2}) and (\ref{basis-kzero3}).
This proves that ${\bf u}$ is semiclassical of class $\mathfrak{s}({\bf u})\leq M+m-1$.
We conclude pursuing with the described procedure,
so that by taking successively $n=0,1,\ldots,m$ in (\ref{basis-kzero1}),
we conclude that the following relations hold:
\begin{equation}\label{basis-kzero4}
D^{m-j} \big(\pi_N {\bf v}\big) = \Phi(\cdot;j)\textbf{u}\quad (j=0,1,\ldots,m) \;.
\end{equation}
In particular, for $j=m$ we obtain (\ref{EqFuc2}), hence ${\bf u}$ and ${\bf v}$
are related by a rational transformation. Next, setting $j=m-1$ in (\ref{basis-kzero4})
we obtain
\begin{equation}\label{basis-kzero5}
D\big(\pi_N {\bf v}\big) = \Phi(\cdot;m-1)\textbf{u}\;.
\end{equation}
Since
$D\big(\Phi(\cdot;m)\pi_N{\bf v}\big)=\Phi(\cdot;m)'\pi_N{\bf v}+\Phi(\cdot;m)D\big(\pi_N{\bf v}\big)$,
we obtain (\ref{EqFuc3}) using (\ref{basis-kzero5}) and (\ref{EqFuc2}).
Thus ${\bf v}$ is semiclassical of class $\mathfrak{s}({\bf v})\leq N+M+2m-2$,
and the theorem is proved.
\end{proof}

\begin{remark}
In the case $m=1$, Theorem \ref{theorem-main-kzero} was partially proved in \cite{MBP1993}.
Note that the functional equation (\ref{EqFuc3}) (for $m=1$) was not given therein.
\end{remark}

\begin{remark}
%The results stated in \cite{JMPN2015} for the continuous  $(M,N)-$coherent pairs of order $(m,k)$ were extended in \cite{JCAM2015} to the setting of discrete OPS. In a similar way, the results proved in this section may be extended
%to the discrete OPS, replacing the derivative operator $D$ by $E:=\triangle_\omega/\omega$ or $D_{q}$, defined by
%$$
%\triangle_\omega f(x):=f(x+\omega)-f(x)\;,\quad D_qf(x):=\frac{f(qx)-f(x)}{qx-x}\quad (f\in\mathcal{P})\;,
%$$
%where $q\neq e^{2ij\pi/n}$ for $0\leq j\leq n-1$ and $n=1,2,\ldots$.
Given complex numbers
$\omega$ and $q$ such that $|q-1|+|\omega|\neq0$,
%\begin{equation}\label{q-notexp}
%q\neq 0,{\rm e}^{2ij\pi/n}\quad(0\leq j\leq n-1\;;\;\;n=1,2,\ldots)\;,
%|q-1|+|\omega|\neq0\;,
%\end{equation}
the operator $D_{q,\omega}:\mathcal{P}\to\mathcal{P}$  considered by Hahn in his influential work \cite{H1949} is defined by 
\begin{equation}\label{def-Dqw}
D_{q,\omega}f(x):=\frac{f(qx+\omega)-f(x)}{(q-1)x+\omega}\quad (f\in\mathcal{P})\;.
\end{equation}
The results and proofs in this section can be repeated with almost no changes in the more general setting of the discrete OPS, replacing the derivative operator $D$ by $D_{q,\omega}$. Actually, the same can be done for discrete OPS on a non-uniform lattice.
\end{remark}

\section{An application}\label{Section-example}

Let $(P_n)_{n\geq0}$ be a monic OPS with respect to a positive Borel measure.
Suppose that $(P_n)_{n\geq0}$ satisfies the differential-difference equation
\begin{equation}\label{Griffin1}
\pi(x)P'_n(x)=b_nP_n(x)+(c_nx+d_n)P_{n-1}(x)\quad(n=0,1,\ldots)\;,
\end{equation}
where $\pi(x)$ is a monic polynomial of degree $1$ and $(b_n)_{n\geq0}$, $(c_n)_{n\geq0}$, and $(d_n)_{n\geq0}$
are sequences of real numbers, with $c_n\neq0$ for each $n=1,2,\ldots$.
We assume $$\pi(x)=x\,.$$
OPS characterized by equation (\ref{Griffin1}) have been studied recently in \cite{G2016}.
Here we give an alternative approach based on the general results presented in the previous section.
$(P_n)_{n\geq0}$ is characterized by a three-term recurrence relation:
\begin{equation}\label{TTRRPn}
xP_n(x)=P_{n+1}(x)+\beta_n P_n(x)+\gamma_nP_{n-1}(x)\quad(n=0,1,\ldots)\;,
\end{equation}
where $(\beta_n)_{n\geq0}$ and $(\gamma_n)_{n\geq1}$
are sequences of real numbers such that $\gamma_n>0$ for each $n\geq1$.
We set $P_{-1}(x)=0$ and $\gamma_0:=0$.
Using (\ref{TTRRPn}), we rewrite (\ref{Griffin1}) as
\begin{equation}\label{piC1}
x\, \frac{P_{n+1}'(x)}{n+1}=P_{n+1}(x)+r_nP_n(x)+s_nP_{n-1}(x)\quad(n=0,1,\ldots)\;,
\end{equation}
where 
$$
r_n:=\frac{c_{n+1}\beta_n+d_{n+1}}{n+1}\;,\quad
s_n:=\frac{c_{n+1}\gamma_n}{n+1}\quad (n=0,1,\ldots)\;.
$$
Notice that $s_n\neq0$ for each $n=1,2,\ldots$.
Comparing (\ref{piC1}) with (\ref{1b}), we have
\begin{equation}\label{mainExamp}
N=M=m=1\;,\;\; k=0\;,\;\; c_{n,n+1}=1\;,\;\; c_{n,n}=r_n\;,\;\; c_{n,n-1}=s_n\;.
\end{equation}
Thus $\big((P_n)_{n\geq0},(P_n)_{n\geq0}\big)$ is a $\pi_1-$coherent pair with index $1$ and order $(1,0)$,
where $\pi_1(x)=x$.
By Theorem \ref{theorem-main-kzero}, the functional ${\bf u}$ with respect to which
$(P_n)_{n\geq0}$ is orthogonal satisfies the relations
\begin{align}
& D\big(\Phi(\cdot;1){\bf u}\big)=\Phi(\cdot;0){\bf u} \label{EqFuc1E1}\\
& x{\bf u}=\Phi(\cdot;1){\bf u} \label{EqFuc2E1}\;. %\\
\end{align}
Since ${\bf u}$ is regular, then (\ref{EqFuc2E1}) implies
\begin{equation}\label{Phi1-1}
\Phi(x;1)=x\;.
\end{equation}
On the other hand, by (\ref{EqPhikzero1})
and using the relations $\beta_n=\langle{\bf u},xP_n^2\rangle/\langle{\bf u},P_{n}^2\rangle$ and
$\gamma_{n+1}=\langle{\bf u},P_{n+1}^2\rangle/\langle{\bf u},P_{n}^2\rangle$ $(n=0,1,\ldots)$, we have
\begin{equation}\label{PhikzeroE0}
\Phi(x;0):=-\frac{r_0}{\gamma_1}P_1(x)-\frac{2s_1}{\gamma_1\gamma_2}P_2(x)\;.
\end{equation}
From (\ref{piC1}) for $n=0,1,2$, and taking into account (\ref{TTRRPn}), we deduce
\begin{equation}\label{eqb0b1b2c1c2}
\begin{aligned}
& r_0=\beta_0\;,\quad
r_1=\mbox{$\frac12$}\,(\beta_0+\beta_1)\;,\quad
r_2=\mbox{$\frac13$}\,(\beta_0+\beta_1+\beta_2)\;,  \\
& s_1=\gamma_1+\mbox{$\frac12$}\,\beta_0(\beta_0-\beta_1)\;, \quad
\beta_0 (s_2-\gamma_2)=(\beta_0\beta_1-\gamma_1)(r_2-\beta_2)\;, \\
& s_2=\mbox{$\frac13$}\,\big(\beta_0^2+\beta_1^2-(\beta_0+\beta_1)\beta_2+2(\gamma_1+\gamma_2)\big)\;.
\end{aligned}
\end{equation}
Therefore, taking into account (\ref{Phi1-1})--(\ref{eqb0b1b2c1c2}) and (\ref{TTRRPn}),
(\ref{EqFuc1E1}) reduces to
\begin{equation}\label{EqFunGriffin}
D\big(x{\bf u}\big)=(-2ax^2+bx+c+1){\bf u} \;,
\end{equation}
where
\begin{align*}
& a:=\frac{s_1}{\gamma_1\gamma_2}=\frac{2\gamma_1+(\beta_0-\beta_1)\beta_0}{2\gamma_1\gamma_2}\;, \\
& b:=\frac{\big(2\gamma_1+(\beta_0-\beta_1)\beta_0\big)(\beta_0+\beta_1)-\beta_0\gamma_2}{\gamma_1\gamma_2}\;, \\
& c:=\frac{\beta_0^2\gamma_2-\big(2\gamma_1+(\beta_0-\beta_1)\beta_0\big)(\beta_0\beta_1-\gamma_1)}{\gamma_1\gamma_2}-1\;.
\end{align*}
Using (\ref{eqb0b1b2c1c2}), and assuming $s_1>0$, we deduce
\begin{equation}\label{eqb0b1b2c1c2-reverted}
\begin{aligned}
& \beta_0=r_0\;,\quad
\beta_1=2r_1-r_0\;,\quad
\gamma_1=s_1-r_0(r_0-r_1)\;,  \\
& \gamma_2=\frac{s_1(3s_2-2s_1)+2r_1\big(s_1(2r_0-r_1)-r_0r_1(r_0-r_1)\big)}{2s_1+r_0r_1}\;.
\end{aligned}
\end{equation}
(Notice that $2s_1+r_0r_1\neq0$; indeed, using $\gamma_1=s_1-r_0(r_0-r_1)$, we have
$2s_1+r_0r_1=\gamma_1+s_1+r_0^2>0$.)
Thus $a$, $b$, and $c$ may be written only in terms of $r_0$, $r_1$, $s_1$, and $s_2$.
Hereafter we impose the (integrability) conditions
\begin{equation}\label{ac}
a>0\;,\quad c>-1\;.
\end{equation}
(Note that the condition $a>0$ is equivalent to $s_1>0$ in equation (\ref{piC1}), or to $c_2>0$ in 
equation equation (\ref{Griffin1}).)
Let $w$ be a solution of 
\begin{equation}\label{EqODEGriffin2}
xw'(x)=(-2ax^2+bx+c)w(x)\;,\quad x\in \mathbb{R}\setminus\{0\}\;.
\end{equation}
Solving this equation imposing (without loss of generality) $w$ to be right-continuous at $x=0$, we find
\begin{equation}\label{ODEw}
w(x)=\left\{
\begin{array}{lcl}
K_1 |x|^c e^{-ax^2+bx} & \mbox{\rm if} & x<0 \;, \\ [0.25em]
K_2 |x|^c e^{-ax^2+bx}  & \mbox{\rm if} & x\geq0 \;,
\end{array}
\right.
\end{equation}
$K_1$ and $K_2$ being real constants.
Requiring, in addition, $K_1$ and $K_2$ to be non-negative and no simultaneously equal to zero,
$w$ becomes a weight function, i.e.,
a non-negative and integrable function which does not vanishes identically and
having finite moments of all orders.
Now, define a functional ${\bf w}$ by
$$
\langle{\bf w},f\rangle:=\kappa\int_\mathbb{R}f(x)w(x)\,{\rm d}x\quad (f\in\mathcal{P})\;,
$$
where $\kappa$ is a normalization constant chosen so that $\langle{\bf w},1\rangle=\langle{\bf u},1\rangle$.
Using (\ref{EqODEGriffin2}) and integration by parts, together with the rules of the distributional calculus,
we show that $D\big(x{\bf w}\big)=(-2ax^2+bx+c+1){\bf w}$ on $\mathcal{P}'$, hence
${\bf w}$ fulfills the same functional equation (\ref{EqFunGriffin}) as ${\bf u}$.
This is equivalent to saying that the sequences of moments $(u_n)_{n\geq0}$ and $(w_n)_{n\geq0}$ of
${\bf u}$ and ${\bf w}$ (defined by $u_n:=\langle{\bf u},x^n\rangle$ and $w_n:=\langle{\bf w},x^n\rangle$)
are solutions of the second order linear difference equation
$$
-2av_{n+2}+(n+b)v_{n+1}+(c+1)v_n=0\quad(n=0,1,\cdots)\;.
$$
Now we show that we may choose $K_1$ and $K_2$ so that ${\bf u}={\bf w}$.
Indeed, since by definition of ${\bf w}$ the condition $u_0=w_0$ holds, we only need to show that
we may choose $K_1$ and $K_2$ so that $u_1=w_1$.
Indeed,
$$
\kappa^{-1}w_1=\int_\mathbb{R}xw(x)\,{\rm d}x
=K_1\int_{-\infty}^0x|x|^{c}e^{-ax^2+bx}\,{\rm d}x+K_2\int_0^{+\infty}x^{c+1}e^{-ax^2+bx}\,{\rm d}x\;,
$$
and making the change of variables $x\mapsto -x$ on the first integral, we obtain
$$
w_1=\kappa \Big(K_2\int_0^{+\infty}x^{c+1}e^{-ax^2+bx}\,{\rm d}x
-K_1\int_0^{+\infty}x^{c+1}e^{-ax^2-bx}\,{\rm d}x\Big)\;.
$$
On the other hand, from $P_1(x)=x-\beta_0$, we have $u_1=\beta_0u_0=r_0w_0$, i.e.,
$$
u_1=\kappa r_0\Big(K_2\int_0^{+\infty}x^{c}e^{-ax^2+bx}\,{\rm d}x
+K_1\int_0^{+\infty}x^{c}e^{-ax^2-bx}\,{\rm d}x\Big)\;.
$$
Therefore, in order to have $u_1=w_1$, we need to impose
$$
r_0=\frac{K_2\int_0^{+\infty}x^{c+1}e^{-ax^2+bx}\,{\rm d}x
-K_1\int_0^{+\infty}x^{c+1}e^{-ax^2-bx}\,{\rm d}x}{K_1\int_0^{+\infty}x^{c}e^{-ax^2-bx}\,{\rm d}x
+K_2\int_0^{+\infty}x^{c}e^{-ax^2+bx}\,{\rm d}x}\;.
$$
Assuming without loss of generality that $K_2>0$, and setting $M=K_1/K_2$,
this is achieved provided that
\begin{equation}\label{defM}
M=\frac{\int_0^{+\infty}x^{c+1}e^{-ax^2+bx}\,{\rm d}x-r_0 \int_0^{+\infty}x^{c}e^{-ax^2+bx}\,{\rm d}x}{\int_0^{+\infty}x^{c+1}e^{-ax^2-bx}\,{\rm d}x+r_0\int_0^{+\infty}x^{c}e^{-ax^2-bx}\,{\rm d}x}\;.
\end{equation}
Thus, up to a positive constant factor, ${\bf u}$ admits the integral representation
$$
\langle{\bf u},f\rangle:=\int_\mathbb{R}f(x)w(x)\,{\rm d}x\quad (f\in\mathcal{P})\;.
$$
We remark that $w$ is a.e. on $\mathbb{R}$ the unique weight function with respect to which $(P_n)_{n\geq0}$
is a monic OPS. This is an immediate consequence of the fact that the moment problem associated to the distribution function with weight $w$ is determined, as we may see easily taking into account Riesz uniqueness criterium (see e.g. \cite[Theorem II-5.2]{F1971}).
Finally, set
\begin{equation}\label{defuMtc}
{\bf u}^{(M,t,c)}:=h_{\sqrt{a}}{\bf u}\;,\quad t:=b/\sqrt{a}\;,
\end{equation}
meaning that
$\langle{\bf u}^{(M,t,c)},x^n\rangle:=\langle{\bf u},\big(\sqrt{a}\,x\big)^n\rangle$
for each $n=0,1,\ldots$.
Note that making the change of variables $x\to x/\sqrt{a}$ in the integrals appearing in (\ref{defM})
we obtain
\begin{equation}\label{defMt}
M=\frac{\int_0^{+\infty}\big(x-\sqrt{a}\,r_0\big)x^{c}e^{-x^2+tx}\,{\rm d}x}{\int_0^{+\infty}\big(x+\sqrt{a}\,r_0\big)x^{c}e^{-x^2-tx}\,{\rm d}x}\;.
\end{equation}
Since ${\bf u}$ fulfils (\ref{EqFunGriffin}) then ${\bf u}^{(M,t,c)}$ satisfies
$$
D\big(x{\bf u}^{(M,t,c)}\big)=(-2x^2+tx+c+1){\bf u}^{(M,t,c)}\;.
$$
Let $(P^{(M,t,c)}_n)_{n\geq0}$ be the monic OPS with respect to ${\bf u}^{(M,t,c)}$. Then (\ref{defuMtc}) implies
\begin{equation}\label{PnMtc}
P_n(x):=\frac{1}{(\sqrt{a}\,)^n}\,P_n^{(M,t,c)}\big(\sqrt{a}\,x\big)\quad(n=0,1,\ldots)\;.
\end{equation}
Moreover, up to a constant factor, ${\bf u}^{(M,t,c)}$ admits the integral representation
$$
\langle{\bf u}^{(M,t,c)},f\rangle:=\int_\mathbb{R}f(x)w^{(M,t,c)}(x)\,{\rm d}x\quad (f\in\mathcal{P})\;,
$$
where
\begin{equation}\label{ODEwtilde}
w^{(M,t,c)}(x):=\frac{a^{c/2}}{K_2}w\Big(\frac{x}{\sqrt{a}}\Big)=\left\{
\begin{array}{lcl}
M |x|^c e^{-x^2+tx} & \mbox{\rm if} & x<0 \;, \\ [0.25em]
|x|^c e^{-x^2+tx}  & \mbox{\rm if} & x\geq0 \;.
\end{array}
\right.
\end{equation}

In conclusion, if $(P_n)_{n\geq0}$ is a monic OPS with respect to a positive-definite linear functional and fulfills (\ref{piC1}), where $(r_n)_{n\geq0}$ and $(s_n)_{n\geq1}$ are sequences of real numbers such that $s_{n}\neq0$ for each $n=1,2,\ldots$, then $P_n$ is given by (\ref{PnMtc}) ---$(P^{(M,t,c)}_n)_{n\geq0}$ being the unique monic OPS with respect to the weight function $w^{(M,t,c)}$ defined by the right-hand side of (\ref{ODEwtilde})---, provided that conditions (\ref{ac}) hold for each choice of the four (real) parameters $r_0$, $r_1$, $s_1$, and $s_2$.

For instance, choosing $r_0=r_1=0$, $s_1=1/2$, and $s_2=1$, we obtain $a=1$, $t=c=0$, and $M=1$, hence $w^{(1,0,0)}(x)=e^{-x^2}$, so that $(P_n)_{n\geq0}$ is the Hermite monic OPS (up to an affine change of the variable).
Finally, we note that (\ref{ODEwtilde}), (\ref{PnMtc}), and (\ref{defMt}) agree, respectively, with (2.27), (2.29), and (2.30) in \cite{G2016}.

\section*{Acknowledgements }
The authors are indebted to Professor J. Petronilho for suggesting this problem, as well as his time for many very helpful discussions that led to the ideas presented. The authors also gratefully acknowledge fruitful discussions with Professor R. \'Alvarez-Nodarse. This work waspartially supported by the Centre for Mathematics of the University of Coimbra--UID/MAT/00324/2019, funded by the Portuguese Government through FCT/MEC and co-funded by the European Regional Development Fund through the Partnership Agreement PT2020.

\end{document}